\documentclass{article}
\usepackage{verbatim}

\usepackage{amsmath}
\usepackage{amssymb}
\usepackage{latexsym}
\usepackage{staves}
\usepackage{marvosym}
\usepackage{eurosym}
\usepackage{phaistos}
\usepackage{verbatim}

\usepackage{color}
\usepackage[T1]{fontenc}
\usepackage[english]{babel}
\usepackage{graphicx}
\usepackage{tikz}
\usetikzlibrary{arrows}
\usetikzlibrary{arrows,decorations.markings}
\usetikzlibrary{mindmap,trees}
\usepackage{ulem}

\newtheorem{Definition}{Definition}[part]
\newtheorem{Proposition}{Proposition}[part]

\newtheorem{Assumption}{Assumption}[part]

\newtheorem{Remark}{Remark}[part]

\makeatletter \@addtoreset{equation}{section}

\@addtoreset{Definition}{section}

\@addtoreset{Theorem}{section}

\@addtoreset{Proposition}{section}

\@addtoreset{Property}{section}

\@addtoreset{Assumption}{section}

\@addtoreset{Corollary}{section}

\@addtoreset{Lemma}{section}

\@addtoreset{Remark}{section}

\@addtoreset{Example}{section}

\def \Sum{\displaystyle\sum}
\def \Prod{\displaystyle\prod}

\def \proof{{\noindent \bf Proof. }}
\def \ep{\hbox{ }\hfill$\Box$}
\def\reff#1{{\rm(\ref{#1})}}

\def\Ac{{\cal A}}

\def\eps{\varepsilon}

\def\Fc{{\cal F}}

\def\eps{\varepsilon}

\def\diag{\mbox{\rm diag}}

\def\no{\noindent}

\def\x{\times}

\def\05{\frac{1}{2}}
\def\-1{^{-1}}
\def\1{{1\hspace{-1mm}{\rm I}}}
\def\={\;=\;}
\def\.{\;.}

\title{  }
\author{ }

\def\be{\begin{eqnarray}}
\def\ee{\end{eqnarray}}

\def\b*{\begin{eqnarray*}}
\def\e*{\end{eqnarray*}}

\def \ep{\hbox{ }\hfill{ ${\cal t}$~\hspace{-5.5mm}~${\cal u}$   } }

\def\E{\mathbb{E}}
\def\F{\mathbb{F}}

\def\R{\mathbb{R}}

\def\P{\mathbb{P}}

\def \F{I\!\!F}

\def \R{I\!\!R}

\def\Ac{{\cal A}}

\def\Fc{{\cal F}}

\def\-1{^{-1}}

\def\0.5{\frac{1}{2}}

\def\x{\times}
\def\no{\noindent}
\def\={\;=\;}
\def\.{\;.}

\def \proof{{\noindent \bf Proof. }}
\def\reff#1{{\rm(\ref{#1})}}
\def\eps{\varepsilon}

\def\diag{\mbox{\rm diag}}

\def\1{{\bf 1}}

\def \ep{\hbox{ }\hfill{ ${\cal t}$~\hspace{-5.5mm}~${\cal u}$   } }
\def \Sum{\displaystyle\sum}
\def \Prod{\displaystyle\prod}

\title{ An Extended Mean Field Game for Storage in Smart Grids }

\author{Cl\'emence Alasseur \footnote{ \texttt{clemence.alasseur@edf.fr} , {\sc edf r\&d} and Finance for energy Market Research Centre ({\sc f}i{\sc me})}\\
	\and
Imen  Ben Taher
	\footnote{ \texttt{imen@ceremade.dauphine.fr}, {\sc ceremade} and Finance for energy Market Research Centre ({\sc f}i{\sc me}) }
	\\
	\and
	Anis Matoussi
	\footnote{ \texttt{anis.matoussi@univ-lemans.fr}, { Institut du Risque et de l'Assurance, Le Mans U}niversit\'e }
	\\
}

\begin{document}

\maketitle

\begin{abstract}
	
	We consider a stylized model for a power network with distributed local power generation and storage. This system is modeled as a network connection of a large number of  nodes, where each node is  characterized by a local electricity consumption, has a local electricity  production (e.g. photovoltaic panels), and manages a local storage device. Depending on its instantaneous consumption and production rate as well as  its storage management decision, each node  may   either buy or sell electricity, impacting the electricity spot price.    The objective at each node is to minimize energy and storage costs by optimally controlling the storage device.  In a non-cooperative game setting, we are led to the analysis of a non-zero sum stochastic game with $N$ players where the interaction  takes place  through the spot price mechanism.  For an infinite number of agents, our model corresponds to an Extended Mean-Field Game.  We are able to compare this solution to the optimal strategy of a central planner and in a linear quadratic setting, we obtain and explicit solution to the Extended Mean-Field Game and we show that it provides an approximate Nash-equilibrium for  $N$-player game.

	\no {\bf Keywords:~ } smart-grid, distributed generation, stochastic renewable generation, optimal storage, stochastic control, extended mean-field games
\end{abstract}

\section{Introduction}

Until the late $90$'s, the power system was characterized by  predictable supply insured by massive {\sl vertically-integrated utilities} which assumed the three major services: generation, transmission and distribution.  Since then,  critical changes have been occurring, and the centralized and vertically-integrated scheme is giving way to a new scheme  where small-scale  distributed generation and storage have an important weight \cite{DelftReport}. Indeed, technological innovation and environmental concerns triggered and are still boosting  the integration of intermittent renewable energy, part of which is provided by relatively small and geographically distributed generation.  The fast growing deployment of decentralized small scale power generation is aided by the  
simultaneous  evolution of local storage technologies and its complementary deployment. This transition calls for in-depth re-engineering of distribution networks at various levels, including tariff structures. A growing literature is interested in distributed storage management and the analysis of its development within the system. In particular, mean field games (MFG) approach has been already used by \cite{Couillet21012} who analyze a system with controlled electrical vehicles and by \cite{Paola2016} with local batteries. These two papers deal with numerical analysis of the  corresponding MFG without providing the existence and uniqueness of the optimal control results.

\paragraph{Our Mean Field Game model for a  power network with distributed storage and generation.}
 The aim of a our paper is to provide a stylized quantitative model for a power system with distributed local  energy generation and storage where   some questions arising in this power grid can be tractably analyzed. This system is modeled as a  network connecting  a large number of  nodes.   Each node has a local electricity consumption,  a local electricity    production (e.g. photovoltaic panels), and manages a local storage device. In our model each node  is characterized by two state variables: the  {\sl local net production} $Q_t $ and  the {\sl battery level} $S_t$, and a control variable: the storage action $\alpha_t$.  At each moment, $Q_t - \alpha_t$ can be either positive or negative ; if positive, respectively negative,  it corresponds to electricity  that the node sells to, resp. buys from, the grid at the spot price.  We consider that  objective of each  node is to minimize its own cost of electricity consumption by controlling the  storage device. As in \cite{Couillet21012} or \cite{Paola2016},   we  assume that the spot price level reflects the instantaneous global consumption, hence, it depends on the strategies of the nodes.  
 In a non-cooperative game setting, we are led to the analysis of a non-zero sum stochastic game with $N$ players and to the search of  Nash equilibria. By making the hypothesis that $N$ goes to infinity, we rely on a Mean Field Game (MFG) approach, more precisely we formulate and solve an  Extended Mean Field Game (EMFG) with common noise.

\paragraph{Literature review for MFG and FBSDE.}

 First we mention that mean field game theory was introduced by the parallel works of Huang and Malhame   \cite{huang2006,huang2007} and of Lasry and Lions \cite{ lasry2006jeux-I,lasry2006jeux-II}, see the notes of  Cardaliaguet \cite{Cardialiaguet-LN13} based on the lectures of  P.-L. Lions  at the Collège of France \cite{Lions-Lecture07}, and the recent the book of Carmona and Delarue  \cite{CarmonaDelarueBook17}. 
 Carmona, Delarue, and Lacker \cite{CarmonaDelarueLacker16} have developed a probabilistic approach based on a stochastic maximum principle for a representative player and use a fixed point argument to find a mean field Nash equilibrium. 
A related but distinct concept is that of mean field type control. In this case, the goal is to assign a strategy to all players at once, such that the resulting crowd behavior is optimal with respect to costs imposed on a central planner. For a comparison of mean field games and mean field type control, see the book of Bensoussan, Frehse, and Yam \cite{BensoussanFrehseYam-book13} (see also \cite{BensoussanFrehseYam17}) as well as the article by Carmona, Delarue, and Lachapelle \cite{CarmonaDelarueLachapelle13}. A key reference is the work of Carmona and Delarue \cite{CarmonaDelarue15}, which characterizes solutions to the mean field type control problem in terms of a stochastic maximum principle for McKean-Vlasov type dynamics (see also \cite{CarmonaDelarueLacker16},  \cite{CarmonaDelarueLacker17}).\\
 Conceptually,  mean field type control  (MFC) is different from the mean field game (MFG),  and although in general an optimal control  on MFC is not an equilibrium strategy on  MFG,  nevertheless  Lasry and Lions in \cite{lasry2007mean} have pointed out that in many cases a mean field Nash equilibrium is also the solution to  an optimal control problem.
 The work of Graber \cite{Graber2016}  also have highlighted this point of view. Motivated by  economic examples, he reformulated  the Nash equilibrium  for  MFG as an optimal control problem. Therefore, he has studied the mean field type control problem associated to the MFG, even though, a priori, he was interested in mean field games.  The present work also follows this point of view.
 
\paragraph{Main contributions.}

A primary   contribution of this paper is that the EMFG approach provides an analytically and numerically tractable setting to assess questions related to the distributed generation and storage. Under proper conditions, the EMFG we associate to this power network game is proven to admit a unique solution which can be characterized though solving an associated Forward Backward Stochastic Differential Equations (FBSDE).  
In the particular case where the cost structure is quadratic and the pricing rule is linear, the FBSDE which characterizes the solution of the EMFG can even be solved explicitly.  This provides a quite tractable and efficient setting to analyze  numerically various questions arising in this power grid. For example, our model gives indications to the question on how decentralized batteries could spread, be managed and how this will impact the spot price depending on the electricity tariff structure. Our model also points out how characteristics of the prosumers' consumptions/productions such as their  seasonal pattern and their volatility change the way they manage a storage. To our knowledge, only the paper by \cite{CardaliaguetLehalle17} also provides explicit solution for an EMFG applied to optimal liquidation of a portfolio. We refer the reader to \cite{carmona2013} for general discussion on the probabilistic approach for MFG. \\ A secondary, yet important finding, is that  our   EMFG can be    profitably compared to a suitable Mean Field Type Control (MFC) problem whose solution  can be interpreted as the optimal strategy of a central planner who coordinates the storage actions at the nodes. Therefore, our model gives clues to an aggregator on how to manage a collection of consumers in a decentralized way.

\paragraph{Structure of the paper.} This paper is organized as follows. In Section \ref{sect: power grid model} we introduce the  stylized model for the power network,  we define the associated $N$-players Nash Game as well as the problem of a central planner who aims to optimally coordinate the storage in the nodes. In Section \ref{sect: mfg approx} we provide the EMFG approximation, characterize the solution of this EMFG and show how it compares to the solution of the MFC problem related to the central planner. In Section \ref{sect: explicit linear quadratic case} we provide and discuss the explicit solution in the  particular case where the cost structure is quadratic and the pricing rule is linear. Finally, in Section \ref{sect:  numerical} a numerical case of study is detailed:  ~our model is applied to the case where the network is composed by two types of agents:~ group 1 of traditional consumers with no local production nor storage, and group 2  of prosumers with local production and storage.  Both the EMFG and central planner strategies are analyzed, compared and commented.

\section{The power grid model}
\label{sect: power grid model}
We consider  a stylized model for a power grid with  distributed local  energy generation and storage. The grid connects $N$ nodes indexed by $i=1, \cdots,N$. Each node is characterized by two state variables:~the    {\sl local net power production}   $Q^i_t$ which represents the local power production  minus the local power consumption at node $i$, and the storage level $S^i_t$ which represents the  total {\sl energy} available in the storage device. We assume that the nodes forming this grid can be partitioned in $\Gamma$ different groups:  the nodes in the  same group $\gamma$ share same  characteristics of local net power production and storage, yet these characteristics vary from one group  to the other.

\no 
 We denote by $N_\gamma$ the number of nodes in group $\gamma$, so that $N = \sum_{\gamma=1}^\Gamma  N_\gamma$, and let $\pi^\gamma = N_\gamma /N$ be the ratio of the population size of region $\gamma$ to the whole population.  We shall  abusively write $i \in \gamma$ to signify that the node $i$ is in region $\gamma$.
 
 The grid also connects a group,   indexed by $0$,  which is characterized by one state variable, its    {\sl local net power production}   $Q^0_t$, and which  does not possess any storage.

\begin{tikzpicture}[
  mindmap,
  every node/.style={concept, execute at begin node=\hskip0pt},
  root concept/.append style={
    concept color=black, fill=white, line width=1ex, text=black, scale=0.7
  },
  text=black, grow cyclic,
  level 1/.append style={scale= 0.9,level distance=4.5cm,sibling angle=140},
  level 2/.append style={level distance=3cm,sibling angle=45}
]

\node[root concept](g) {Power Grid} 
child[concept color=black] { node[fill = white] {Rest of the world $Q^0$}
}
child[concept color=blue] { node(r1)[fill=white, scale =0.8] {Region 1}
child { node[fill=white](r11) {node 1 $S^{1}$ $|$ $Q^{1}$ } }
child { node[fill=white](r12) {node 2 $S^{2}$ $|$ $Q^{2}$ } }
child { node[fill=white](r13) {node  i  $S^{i}$ $|$ $Q^{i}$ } }
child { node[fill=white](r14) {node  j $S^{j}$ $|$ $Q^{j} $ } }
child { node[fill=white](r15) {node k $S^{k}$ $|$ $Q^{k}$ } }
child { node[fill=white](r16) {node $N_1$ $S^{N_1}$ $|$ $Q^{N_1}$ } }
}
child[concept color=orange] { node(r2)[fill=white, scale =0.8] {Region 2}
child { node[fill=white](r21) {node 1 $S^{N_1+1}$ $|$ $Q^{N_1+1}$ } }
child { node[fill=white](r22) {node 2 $S^{N_1+2}$ $|$ $Q^{N_1+2}$} }
child { node[fill=white](r23) {node i $S^{N_1+i}$ $|$ $Q^{N_1+i}$}}
child { node[fill=white](r24) {node  $N_2$ $S^{2,N_1+N_2}$ $|$ $Q^{2,N_1+N_2}$ } }
};
 \draw[<->] (g) -- (r1);
 \draw[<->] (r1) -- (r11);
 \draw[<->] (r1) -- (r12); 
 \draw[<->] (r1) -- (r13); 
 \draw[<->] (r1) -- (r14); 
 \draw[<->] (r1) -- (r15); 
 \draw[<->] (r1) -- (r16);   
 \draw[<->] (g) -- (r2);
 \draw[<->] (r2) -- (r21);
 \draw[<->] (r2) -- (r22); 
 \draw[<->] (r2) -- (r23); 
 \draw[<->] (r2) -- (r24); 
\end{tikzpicture}

\begin{Remark}[Partitioning of the nodes]{\rm
 Such a partitioning of the nodes   is relevant  for the modelling and analysis of various situations. For instance,  in Section \ref{sect: numerical} we consider a grid with two types of agents, group 1 consists of traditional consumers with no local production nor storage, and group 2 consists of prosumers with local production and storage.  We may also consider a grid with $\Gamma$ different geographical regions, each region beeing characterized by a specific mode of local power production driven by example by specific meteorological conditions etc.}
\end{Remark}

\bigskip

 In order to model the dynamics of the state variables, we consider  a complete probability space $(\Omega, \Fc, \P)$ on which are defined independent Brownian motions $B^0, B^1, \cdots, B^N$.
  We  consider $N$ independent identically distributed (i.i.d.)  random variables $x_0^i=(s^i_0,q^i_0)$ which are independent of $B^0$ and $B^i$. We denote by $\F=\{\Fc_t\}$ the filtration defined 
  by $\Fc_t = \sigma ((s^i_0,  q^0_0, q^i_0), B^0_s, B^i_s, i=1, \cdots N, ~s \le t\}$, and  by  $\F^0=\{ \Fc^0_t\}$ the filtration generated by $B^0$ i.e.  $\Fc^0_t = \sigma( B^0_s, ~s \le t\}$. We denote by $\Ac$ the set of $\F$-adapted real-valued processes $a=\{a_t\}$ such that
  	$
	\E\left[ \int_0^T |a_u|^2 du \right]< \infty
	$.

\bigskip
\no We assume that at node $i$ in the region $\gamma$, the battery level is controlled through a storage action  $\alpha^{\gamma,i} \in \Ac$ according  to
	\b*
	S^{i}_t&=& s^{i}_0 + \int_0^t \alpha^{i}_s ds, 
	\e*
and that, if the node $i$ is in the region $\gamma$, then
the net power production is given by
	\b*
	dQ^i_t&=&\mu^\gamma(t,Q^i_t) dt  + \sigma^\gamma(t,Q^i_t) dB^i_t + \sigma^{\gamma 0}(t,Q^i_t) dB^0_t\;,~~Q^i_0 = q^i_0.
	\e*
The  {\bf net injection} of the node $i$ is
	\b*
	Q^i_t - \alpha^i_t.
	\e*	
It can be either positive or negative. If positive then  it corresponds to electricity being sold from the node $i$ to the grid ;  if negative, then  it corresponds to electricity being  bought by the node $i$ from  the grid.
		
\no The net injection of the rest of the world is given by 
\b*
dQ^0_t&=&\mu^0(t,Q^0_t) dt  + \sigma^{0}(t,Q^0_t) dB^0_t\;,~~Q^0_0 = q^0_0.
\e*
In our model $B^0_t$  represents a common signal which affects  the energy demand of the whole grid.  Then  for each $i$,   $\sigma^{\gamma 0}:~\R \rightarrow \R$ is a  given function which allows to model how  the  node $i$ of region $\gamma$ is  affected by the common signal $B^0_t$. We assume that the rest of the world is only affected by this common signal $B^0_t$.

	\begin{Remark}[Constraints on the storage]{\rm
			In our model we do not enforce constraints on the storage level nor on the injection/withdrawal rates. Indeed, we give priority to finding explicit solutions to our problem in order to analyse the qualitative behavior of the system.  In the numerical examples we considered  we were able to obtain reasonable interpretations and results despite this limitation on the modeling of the storages.
		}
	\end{Remark}
  
\subsection{Electricity spot price}
 We make the assumption that the electricity price per Watt-hour depends on the instantaneous demand. When the strategy   $\alpha= (\alpha^1, \cdots,\alpha^N) \in    \Ac^N$ is implemented the spot price is given by
	\b*
	P^{N,\alpha}_t&=& p\left(-  Q^0_t -  \sum_{i=1}^N \eta (Q^i_t - \alpha^i_t)\right),
	\e*
where $p(\cdot)$ is the exogeneous inverse demand function for electricity, and $\eta$ is a scaling parameter which weights the contribution of each individual node $i$ to the whole system. 
 We model a grid with a large number of  `small' nodes $i$, hence we shall be considering  the limit as  $N \rightarrow + \infty$ and $\eta \rightarrow 0$. Here we assume that 
  	\b*
	  \eta &=& 1/N
	\e*
Hence the spot price  depends on   the  averaged net injections $\frac{1}{N} \sum_{i=1}^N (Q^i_t - \alpha^i_t)$ 
	 \b*
	P^{N,\alpha}_t&=& p\left( - Q^0_t -  \sum_{i=1}^N  \frac{1}{N}(Q^i_t - \alpha^i_t)\right).\e*

\begin{Assumption}
The function $p(\cdot)$ is assumed to be strictly increasing.
\end{Assumption}

\begin{Remark}{\rm  The fact that the spot price depends on the  averaged net injections $\frac{1}{N} \sum_{i=1}^N (Q^i_t - \alpha^i_t)$ is the rationale for our Extended Mean Field Game (EMFG) approximation developed in Section \ref{sect: mfg approx}.  Recall that  $\pi^\gamma = N_\gamma / N$ and notice that the electricity price can be expressed as
	\begin{equation}
	\label{price:MFG}
	 P^{N,\alpha}_t=
	p\left( -  Q^0_t -  \sum_{\gamma=1}^\Gamma   \pi^\gamma  \, \sum_{i \in \gamma}\frac{1}{N_\gamma}(Q^i_t - \alpha^i_t)\right).
     \end{equation}
At a ``macroscopic level', each region $\gamma$ influences the price through the {\bf average net injection} $\sum_{i \in \gamma}\frac{1}{N_\gamma}(Q^i_t - \alpha^i_t)$ modulated by the ratio $\pi^\gamma = N_\gamma/N$.
 	}
\end{Remark}

\subsection{Cost functions}  
We consider  a finite time horizon $T>0$. When the control action   $\alpha = (\alpha^1, \cdots, \alpha^N)$ is implemented, 
the  cost incurred at the node $i$ in the region $\gamma=1, \cdots, N$ breaks down into three components :~a volumetric charge,  a demand charge, and a storage cost. The first two components correspond to the electricity bill. Indeed, the consumer's bill is commonly the sum of this two components: one proportional to the energy consumed (the volumetric charge) and one linked to the maximum  power achieved (the demand charge). For example, in France, the consumer subscribes to a maximum power level, meaning that its instantaneous consumption is limited to this level physically by its meter and its demand charge component increases with the level subscribed. In other countries, such as in some states of the US, the demand charge could be based on the highest 15-minute average usage recorded on the demand meter within a given month. The third component of the cost corresponds to the costs of the storage (purchase, maintenance, wear). 
\b* 
	J^{i,\gamma,N}(\alpha)&=&\underbrace{\E\left[
	\int_0^T P^{N,\alpha}_t.\left( \alpha^i_t - Q^i_t\right)dt
		\right]  }_{\mbox{volumetric charge}}+ 
		\underbrace{\E\left[
	\int_0^T  L^\gamma_T(Q^{i}_t, \alpha^i_t)dt\right]}_{\mbox{demand charge}}
		\\&& \hspace{23mm}
		+~~\underbrace{\E
	\left[\int_0^T  L_S(S^{i, \alpha^i}_t, \alpha^i_t) 
		dt + g(S^{i, \alpha^i}_T)\right]}_{\mbox{storage cost}}.
	\e*
where $L_T^\gamma, L_S:~\R \x \R \rightarrow \R$ , and $g:~\R \rightarrow \R$ are  continuous functions.

The term $P^{N,\alpha}_t .\left( \alpha^i_t - Q^i_t\right)$ represents  the current volumetric cost  (or profit) of electricity consumed  (or produced)  at the spot price $P^{N,\alpha}_t$.  The term $L^\gamma_T(Q^{i}_t, \alpha^i_t)$ is linked to the maximum instantaneous power. This demand charge component is designed to reflect the fact that electricity system costs are closely related to power the system requires in peak hours: production installed capacities and network are indeed designed to satisfy highest level of peak demand. 
 The term $L_S(S^{i, \alpha^i}_t, \alpha^i_t)$  represents the current storage cost and is assumed to be identical in all the regions $\gamma$.   The  terminal cost $g(S^{i, \alpha^i}_T)$ typically guarantees a minimal level of storage at the end of the period.

\no Finally, the region $0$/ rest of the world incurs  only  energy and transmission costs 
	\be\label{eq: current energy cost r}
	J^{0,N}(\alpha)&=&\underbrace{\E\left[
	\int_0^T  - P^{N,\alpha}_t . Q^0_t dt
		\right]  }_{\mbox{volumetric charge}}+ 
		\underbrace{\E\left[
	\int_0^T  L^0_T(Q^{0}_t, 0)dt
		\right]  }_{\mbox{demand charge}}
	\ee

\begin{Assumption}\label{asmp: cost convex}
	The current cost $(s,q, \alpha) \mapsto L^\gamma_T(q,\alpha) + L_S(s,\alpha)$ is strictly convex with respect to $(s,\alpha)$.
	The terminal cost $s \mapsto g(s)$ is strictly convex with respect to $s$.
\end{Assumption}
	
\begin{Assumption}\label{asmp: cost coercive}
	There exists some constant $C>0$ such that
		\b*
		 | L_T^\gamma(q,a)  | +   |L_S(s,a) | +  |g(s)  | &\le C\left( |q|^2 + |s|^2+ |a|^2\right)	
		\e* 
	\end{Assumption}
\begin{Assumption}\label{asmp: cost differentiable}
	The functions $L^\gamma_T$, $L_S$  and $g$ are  continuously differentiable and their derivatives are a Lipschitz continuous functions.	\end{Assumption}\begin{Assumption}\label{asmp: diffusion coefficients}	
	The coefficients $ \mu^0 (.,.) $ and   $  \sigma^0(.,.)$ (respectively $ \mu^\gamma (.,.)$ and $  \sigma^\gamma (.,.)$)   are  Lipschitz continuous functions and with linear growth in the state variable.
	
\end{Assumption}

\begin{Remark}[On the quadratic costs hypothesis]
{\rm
Though characterization results of EMFG equilibria for more general cost functions exist in the literature, see e.g. \cite{carmona2013}, very few are the cases where a tractable analysis can be worked out, especially under the presence of common noise. In  the case  of  quadratic cost functions, considered in Section \ref{sect: explicit linear quadratic case}
  we are able to provide a quasi-explicit solution which allows to perform   an easy-to-implement numerical analysis for the system. We claim that, even under the quadratic cost assumption, our model can to some extent accommodate  for some relevant cases of study. }
\end{Remark}

\subsection{Optimality criteria}
\paragraph{Non-cooperative game point of view}
The aim of each node $i$ is to minimize the  cost of electricity consumption by controlling the size and the management of the storage device.  In a non-cooperative game setting, 
we are led to the analysis of a non-zero sum stochastic game with $N$ players and  to   the search of  Nash-equilibria:
\begin{Definition}[Nash equilibrium for the $N$-players game] We say that \\ $\alpha^\star=(\alpha^{\star,1}, \cdots, \alpha^{\star,N}) $ belongs to $\Ac^N$ is a  Nash-equilibrium if for each $(i,\gamma)$, for any $u \in \Ac$:~
	$$
	J^{i,\gamma, N}(\alpha^{\star,1}, \cdots, \alpha^{\star,i-1},u, \alpha^{\star,i+1} , \cdots,\alpha^{\star,N})
	\ge
	J^{i,\gamma, N}(\alpha^{\star,1}, \cdots,  \alpha^{\star, N}).$$
	
\end{Definition} 
\begin{Definition}[$\eps$-Nash equilibrium for the $N$-players game]
Let $\eps >0$. We say that $\alpha^\star=(\alpha^{\star,1}, \cdots, \alpha^{\star,N}) \in \Ac^N$ is a  $\eps$-Nash-equilibrium if for each $(i,\gamma)$, for any $u \in \Ac$:~
	$$
	J^{i,\gamma, N}(\alpha^{\star,1}, \cdots, \alpha^{\star,i-1},u, \alpha^{\star,i+1} , \cdots,\alpha^{\star,N})
	\ge
	J^{i,\gamma,N}(\alpha^{\star,1}, \cdots,  \alpha^{\star, N}) - \eps
	.$$
\end{Definition} 

\paragraph{Central Planner point of view} We should also consider the power grid model from the perspective of a central planner whose aim is to dictate a storage rule:~$\alpha= ( \alpha^1, \cdots,\alpha^N)$ 
in order to minimize the {\it egalitarian} cost function between the nodes and the rest of the world
	\b*
	J^{\rm C, N}\left(\alpha\right)=
	J^{0,N}(\alpha)+ \sum_{i=1}^{N}	\frac{1}{N}  J^{i,\gamma,N}(\alpha).
	\e*
where $1/N$ is the scaling parameter which weights the contribution of each individual node to the system. The cost function $J^{C,N}(\alpha)$ can also be written as	
	\b*
	J^{\rm C, N}\left(\alpha\right)=  J^{0,N}(\alpha)+  \sum_{\gamma =1}^\Gamma \pi^\gamma \sum_{i=1}^{N_\gamma}	\frac{1}{N_\gamma} J^{i,\gamma,N}(\alpha).
	\e*
\begin{Definition}[Optimal coordinated plan]
We say that $ \hat \alpha=(\hat \alpha^{1}, \cdots, \hat \alpha^{N}) \in \Ac^N$ is an optimal coordinated plan if:~
	$
	\hat \alpha = {\rm argmin}_{\alpha \in \Ac^N} \; J^{C, N, \eta}(\alpha)
	$.
\end{Definition}

\section{An Extended  Mean Field Game approximation}
\label{sect: mfg approx}

In this section we consider  on the filtered probability space $(\Omega, \Fc, \P, \F)$,  $\Gamma$ Brownian motions $B^\gamma, \gamma =1, \cdots, \Gamma$ which are mutually independent and independent from the Brownian filtration $\F^0$.

\no We shall use the following notation. If $\xi = \{\xi_t \}$ is an $\F$-adapted process, then $\bar \xi = \{ \bar \xi_t\}$ denotes the process defined by :~$\bar \xi_t := \E[\xi_t | \Fc^0_t]$.
\bigskip

Let $x_0=(s_0, q_0) = \left( x_0^{\gamma}=(s_0^{\gamma}, q_0^{\gamma})\right)_{ 1\le {\gamma} \le \Gamma}$ be a random vector which is independent from $\F^0$. Let $Q^0$ and $Q^\gamma$ be the processes defined by
	\be
	\label{eq: q}
	Q^{\gamma}&=& q_0^{\gamma}  + \int_0^t  
	\mu^{\gamma}(u,Q^{\gamma}) du  +\int_0^t   \sigma^{\gamma}(u,Q^{\gamma}) dB^{\gamma}_u +\int_0^t   \sigma^{{\gamma},0}(u,Q^{\gamma}) dB^0_u\\
	\label{eq: q_r}
	Q^0_t&=&q_0^{0}+  \int_0^t \mu^r(u,Q^0_t) du +  \int_0^t \sigma^{0}(u,Q^0_u) dB^0_u\.
	\ee	
If  $\bar \nu =  (\bar \nu^{1}, \cdots, \bar \nu^{\Gamma})$
 is an $\F^0$-adapted   $\R^\Gamma$-valued process, we denote
 	\begin{equation}
	\label{price:MFC}
	P^{\bar \nu}_t=	p\left(-  Q^0_t   - \sum_{{\gamma} \in \Gamma}\pi^{\gamma}\left(\E[ Q^{{\gamma}}_t | \Fc^0_t] - \bar \nu^{{\gamma}}_t \right)\right).
        \end{equation}
 Our approach is based on the following idea: for a large number of nodes $N$, we approximate the electricity price defined by \eqref{price:MFG} by the expression given by \eqref{price:MFC}  and we fix  the mean field control represented by the process $ \bar \nu $. We can then characterize the solution and control process $\alpha$ and we achieve a Nash equilibrium if the conditional expectation of the control process  $\alpha$ is equal to $ \bar \nu $.  Moreover, we  also fix  the conditional expectation of the net power production, but since  this is uncontrolled, we can directly set it to be  $ \E[ Q^{{\gamma}}_t | \Fc^0_t]$.

We now consider two types of  cost functions,  for any control process $\alpha=(\alpha^1, \cdots, \alpha^\Gamma)$ and for each $\gamma =1, \cdots, \Gamma$, 
	\be
	 J^{\gamma}_{x_0}(\alpha^\gamma, \bar{ \nu})=
	\E\int_0^T \left[ P^{\bar \nu}_t (\alpha^\gamma_t - Q^\gamma_t) + 
	L^\gamma_T(Q^\gamma_t, \alpha^\gamma_t) + L_S(S^\gamma_t, \alpha^\gamma_t)\right]dt + \E\left[ g(S^\gamma_t)\right]	\\
	\mbox{and}~~
	J^{C}_{x_0}(\alpha)=
	\E \int_0^T \left[- 
	P^{\bar{\alpha}}_t Q^0_t +L^0_T(Q^0_t, 0 ) \right]dt 
	+ \Sum_{\gamma=1}^\Gamma
	  \pi^\gamma J_{x_0}^\gamma (\alpha^\gamma, \bar \alpha_t)\\
	  \mbox{where}~~\label{eq: fsde}
	S^{{\gamma}}_t = s_0^{\gamma} + \int_0^t \alpha^\gamma_u du.
	\ee 

\begin{Definition}[Mean field Nash equilibrium]
Let $x_0=(s_0,q_0)$ be  a random vector independent from $\F^0$. We say that $\alpha^\star = \{\alpha^{\gamma,\star},1\le {\gamma} \le \Gamma\}$ is a mean field Nash equilibrium   if,   for each $\gamma$, $ \alpha^{\gamma,\star}$ minimizes the function
$\alpha^\gamma \mapsto J^{\gamma}_{x_0}(\alpha^\gamma, \{ \E[\alpha^\star_t | \Fc^0_t]\})$.
\end{Definition}

\begin{Definition}[Mean field optimal control]
Let $x_0=(s_0,q_0)$ be  a random vector independent from $\F^0$. We say that $\hat \alpha= \{\hat \alpha^{\gamma},1\le {\gamma} \le \Gamma\}$ is a mean field optimal control   if,    $ \hat \alpha$ minimizes the function
$\alpha \mapsto J^{C}_{x_0}(\alpha)$.
\end{Definition}

\begin{Proposition}[{Characterization of  mean field Nash equilibria}]
\label{prop: charac mfg 1}
Let $\bar \nu$ be a given $\F^0$-adapted   $\R^\Gamma$-valued process, and $x_0=(s_0, q_0) = \{ x_0^{\gamma}=(s_0^{\gamma}, q_0^{\gamma}), 1\le \gamma \le \Gamma\}$ be a random vector which is independent form $\F^0$.  
Then there exists a unique  control $\alpha^\star = (\alpha^{1,\star}, \cdots, \alpha^{\Gamma, \star}) = \alpha^\star(\bar \nu, x_0)$ such that:~for each $\gamma$, $\alpha^{\gamma,\star}$ minimizes the function   $\alpha^\gamma \mapsto J^{\gamma}_{x_0}(\alpha^\gamma, \bar \nu)$.
Moreover, if $(S^{\gamma,\star}, Q^\gamma)$  is the state process  corresponding to the initial data condition $x^\gamma_0$, to  the control  $\alpha^{\gamma,\star}$,  and to the dynamic \reff{eq: fsde}-\reff{eq: q}, then there exists a unique adapted solution $(Y^{\gamma, \star}, Z^{0,\gamma, \star}, Z^{\gamma, \star})$ of the BSDE
	\be\label{eq: bsde}
	\left\{
	\begin{array}{l l l}
	dY^{\gamma, \star}_t&=& -\partial _s L_S(S^{\gamma,\star}_t, \alpha^{\gamma,\star}_t)dt +   Z^{0,\gamma, \star}_t dB^0_t + Z^{\gamma,\star}_t dB^\gamma_t\\
	Y^{\gamma, \star}_T&=& \partial_s g(S^{\gamma,\star}_T)
	\end{array}
	\right.
	\ee 
satisfying the coupling condition
	\be\label{eq: coupling cnd}
	 0&=&Y^{\gamma,\star}_t + 
	P^{\bar \nu}_t +
	\partial_\alpha L^\gamma_T(Q^\gamma_t,\alpha^{\gamma,\star}_t) + \partial_\alpha L_S(S^{\gamma,\star}_t,\alpha^{\gamma,\star}_t).
		\ee	
Conversely, assume that there exists $(\alpha^{\gamma,\star}, S^{\gamma,\star}, Y^{\gamma, \star}, Z^{0,\gamma,\star},Z^{\gamma,\star})$ which satisfy the coupling condition \reff{eq: coupling cnd}  as well as the FBSDE  \reff{eq: fsde}-\reff{eq: q}-\reff{eq: bsde}, then $ \alpha^{\gamma,\star}$ is the optimal control minimizing $J^{\gamma}_{x_0}(\alpha^\gamma, \bar \nu)$ and $S^{\gamma,\star}$ is the optimal trajectory. If in addition:
 	\be
	\E\left[  \alpha^{\gamma,\star}_t| \Fc^0_t\right] = \bar \nu^{\gamma,0}_t, ~~\forall \gamma =1, \cdots, \Gamma,
	\ee	
then $ \alpha^\star$ is a mean field Nash 	 equilibrium.
\end{Proposition}
\proof The proof  is based on the classical Pontryagin's maximum principle where the characterization of the Mean field Nash equilibrium is given by the associated McKean-Vlasov FBSDEs  \reff{eq: fsde}-\reff{eq: q}-\reff{eq: bsde}.   
Fix some $\gamma \in \{ 1, \cdots, \Gamma\}$. Assumptions \ref{asmp: cost convex}, \ref{asmp: cost coercive}  and \ref{asmp: cost differentiable} ensure that the function 
	\b*
	\alpha^\gamma \in \Ac \mapsto J^{\gamma}_{x_0}( \alpha^\gamma, \bar \nu)\e*
	 is a strictly convex coercive function and Gateaux-differentiable.  
The Gateaux derivative of $J := J^{\gamma}_{x_0}( \cdot, \bar \nu)$ is
	\b*
	d_\beta J^{\rm}(\alpha^\gamma)&=&
	\E\left[ 
	\int_0^T   \left\{ P^{\bar \nu}_u +
	\partial_\alpha L^\gamma_T(Q^\gamma_u,\alpha^\gamma_u) + \partial_\alpha L_S(S^\gamma_u,\alpha^\gamma_u) \right\} \beta_u du \right]\\
	&&  \hspace{2.1cm}
	+\E\left[\int_0^T \partial_s L(S^{\gamma}_u, \alpha^\gamma_u) \tilde S^\beta_u  du 
	\,+ \tilde S^\beta_T\partial_s g(S^{\gamma}_T)
	\right],
	\e*
where $\tilde S^\beta_u$ is the process defined by
	\b*
	d \tilde S^\beta_u &=&\beta_u du,~~\tilde S^\beta_0=0\;.
	\e*
Hence, 	there exists a unique optimal control $\alpha^{\gamma,\star} =\alpha^{\gamma,\star}(\bar \nu, x_0)$ which satisfies the Euler  optimality condition
	\be\label{eq: optimality cnd} \nonumber
	0&=&\E\left[ 
	\int_0^T   \left\{ P^{\bar \nu}_u+
	\partial_\alpha L^\gamma_T(Q^\gamma_u,\alpha^\gamma_u) + \partial_\alpha L_S(S^\gamma_u,\alpha^\gamma_u) \right\} \beta_u du \right]\\
	&&  \hspace{2.1cm}
	+\E\left[\int_0^T \partial_s L(S^{\gamma}_u, \alpha^\gamma_u) \tilde S^\beta_u  du 
	\,+ \tilde S^\beta_T\partial_s g(S^{\gamma}_T)
	\right]
	\ee
Let $S^{\gamma,\star}$ be the associated optimal trajectory, and let $( Y^{\gamma,\star},  Z^{0,\gamma,\star},  Z^{\gamma,\star})$ be the solution to the BDSE \reff{eq: bsde}, then by It\^o Lemma, for each $\beta$
	\be
	\E\left[ \tilde S^\beta_T Y^{\gamma,\star}_T\right]&=&\E\left[
	\int_0^T \left( Y^{\gamma,\star}_t \beta_t - 
	\partial_sL_S(S^{\gamma,\star}_t, \alpha^\star_t)\tilde S^\beta_t
	 \right) dt 
	\right].
	\ee
Taking into account the terminal condition $Y^\star_T = \partial_s g(S^{\gamma,\star}_T)$ and the optimality condition \reff{eq: optimality cnd}, the previous equation leads to 
	\begin{equation}
	\E\left[ 
	\int_0^T 
	\left( Y^{\gamma,\star}_u + 
	P^{\bar \nu}_u +
	\partial_\alpha L^\gamma_T(Q^\gamma_u,\alpha^{\gamma,\star}_u) + \partial_\alpha L_S(S^{\gamma,\star}_u,\alpha^{\gamma,\star}_u)
	\right) \beta_udu
	\right]=0.
	\end{equation}
Since $\beta$ is arbitrary we conclude to the coupling condition \reff{eq: coupling cnd}.

\no Conversely, if $(\alpha^{\gamma,\star},S^{\gamma,\star}, Y^{\gamma, \star}, Z^{0,\gamma,\star},Z^{\gamma,\star})$ satisfies the coupling condition \reff{eq: coupling cnd} and  the FBSDE system  \reff{eq: fsde}-\reff{eq: q}-\reff{eq: bsde}, then we verify that the gateau derivative of $J^{\gamma, \rm MFG}_{x_0}( \cdot, \bar \nu)$ at $\alpha^{\gamma,\star}$ is equal to zero and we conclude by the strict convexity of $J^{\gamma, \rm MFG}_{x_0}( \cdot, \bar \nu)$ to the desired result. 		
\ep

\bigskip


\begin{Proposition}[{Characterization of  mean field optimal controls}]\label{prop: charac mfc} 
Assume that $\hat \alpha=(\hat \alpha^1, \cdots, \hat \alpha^\Gamma)$  minimizes the functional $J^{C}_{x_0}(\alpha)$, and denote by  $\hat S = (\hat S^\1, \cdots, \hat S^\Gamma)$ the corresponding controlled trajectory.  Then there exists a unique adapted solution $(\hat Y = (\hat Y^1, \cdots \hat Y^\Gamma_t),\hat Z = (\hat Z^1, \cdots, \hat Z^\Gamma),\hat Z^0 = (\hat Z^{0,1}, \cdots, \hat Z^{0,\Gamma}))$ of the BSDE
	\be\label{eq: bsde mfc}
	\left\{
	\begin{array}{l l l}
	d\hat Y^\gamma_t&=& - \partial_s L_S(\hat S^\gamma_t, \hat \alpha^\gamma_t) dt  +  \hat Z^{0,\gamma}_t dB^0_t + \hat Z^\gamma_t dB^\gamma_t\\
	\hat Y^\gamma_T&=&\partial_sg(\hat S^\gamma_T)
	\end{array}
	\right.
	\ee 
satisfying the coupling condition:~for all $\gamma= 1, \cdots, \Gamma$
	\be
	\nonumber
	0&=& \hat Y^\gamma_t 
	+\partial_\alpha L^\gamma_T(Q^\gamma_t, \hat \alpha^\gamma_t) 
	+ \partial_\alpha L_S(\hat S_t, \hat \alpha^\gamma_t) + P^{\bar{\hat  \alpha}}_t  \\
	\label{eq: coupling mfc}
	&& \hspace{7mm}
	 -  p'\left(- Q^{0}_t - 
	 \Pi_{\Gamma} \cdot \left(\bar Q_t - \bar {\hat \alpha}_t \right)
	  \right)
	\left( - Q^{0}_t - \Pi_{\Gamma} \cdot \left(\bar Q_t - \bar{\hat  \alpha}_t \right)  \right) 
	\ee

	with $\bar{\hat  \alpha}_t = \E[\hat \alpha_t | \Fc_t^0]$ and $\Pi_{\Gamma}=(\pi_1, \cdots, \pi_\Gamma)^T$.

\no Conversely, suppose $(\hat S, \hat \alpha, \hat Y, \hat Z^0, \hat Z)$ is an adapted solution to the forward backward system \reff{eq: fsde}-\reff{eq: bsde mfc}, with the coupling condition \reff{eq: coupling mfc}, then $\hat \alpha$ is the optimal control minimizing $J^{\rm MFC}_{x_0}(\alpha)$ and $\hat S$ is the optimal trajectory. 
\end{Proposition}
\proof We only prove the necessary condition of Pontryagin's maximum principle for optimality. The sufficient condition could be proven exactly as it is done in Proposition 3.1.
 Assumption \ref{asmp: cost differentiable} insures that the cost function   $\alpha \in \Ac \mapsto J^{\rm C}_{x_0} (\alpha)$ is Gâteaux differentiable.
   with Gateaux derivative given by 
	\b*
	d_\beta J^{C}_{x_0}(\alpha)&=& \Sum_{\gamma } \pi^\gamma \E\left[ \partial_s g(S^\gamma_T)  \tilde S^{\beta^\gamma}_T+ \int_0^T  \partial_s L_S(S^\gamma_u, \alpha^\gamma_u)S^{\beta^\gamma}_u du \right]\\
	 && + \Sum_\gamma \pi^\gamma\E\left[
	\int_0^T \left\{
	 P^{\bar \alpha}_u  + \partial_\alpha L^\gamma_T(Q^\gamma,\alpha^\gamma_t) + 
	 \partial_\alpha L_S(S^\gamma, \alpha^\gamma)
	   \right\}
	    \beta^\gamma_u  du \right]\\
	&& - \sum_{\gamma} \pi^\gamma \E\left[ \int_0^T 
	 p'\left(-  Q^{0}_u - 
	 \Pi_{\Gamma} \cdot \left(\bar Q_u - \bar \alpha_u \right)
	  \right)
	\left( - Q^{0}_u - \Pi_{\Gamma} \cdot \left(\bar Q_u - \bar \alpha_u \right) 
	 \right\}  \beta^\gamma_udu
		\right],
	\e*

		
where $\tilde S^{\beta^\gamma}_u$ is the process defined by
	\b*
	d \tilde S^{\beta^\gamma}_u &=&\beta^\gamma_u du,~~\tilde S^{\beta^\gamma}_0=0\;.
	\e*
Hence the optimal control $\hat \alpha$ satisfies the Euler optimality condition:~for all $\beta = (\beta^1, \cdots, \beta^\Gamma) $
	\b*
	0&=& \Sum_{\gamma } \pi^\gamma \E\left[ \partial_s g(S^\gamma_T)  \tilde S^{\beta^\gamma}_T+ \int_0^T  \partial_s L_S(S^\gamma_u, \alpha^\gamma_u)S^{\beta^\gamma}_u du \right]\\
	 && + \Sum_\gamma \pi^\gamma\E\left[
	\int_0^T \left\{
	 P^{\bar \alpha}_u  + \partial_\alpha L^\gamma_T(Q^\gamma_u,\alpha^\gamma_u) + 
	 \partial_\alpha L_S(S^\gamma_u, \alpha^\gamma_u)
	   \right\}
	    \beta^\gamma_u  du \right]\\
	&& - \sum_{\gamma} \pi^\gamma \E\left[ \int_0^T 
	\left\{ p'\left(-  Q^{0}_u - 
	 \Pi_{\Gamma} \cdot \left(\bar Q^0_u - \bar \alpha^0_u \right)
	  \right)
	\left(  -  Q^{0}_u -  \Pi_{\Gamma} \cdot (\bar Q_u - \bar \alpha_u \right) 
	 \right\}  \beta^\gamma_udu
		\right],
	\e*  
Now, let $(\hat Y, \hat Z, \hat Z^0)$ be the unique solution to the BSDE \reff{eq: bsde mfc}, and let $\hat S$ be the state process associated to the optimal control $\hat \alpha$, applying It\^o formula, we obtain 	
	\b*
	\sum_{\gamma}  \pi^\gamma \E\left[ \hat Y^\gamma_T \tilde S^{\beta^\gamma}_T\right]&=&
	\sum_{\gamma} \pi^\gamma\E\left[  
	\int_0^T 
	\left\{- \partial_s L_S(\hat S_u, \hat \alpha_u) + \beta^\gamma_u \hat Y^\gamma_u  \right\}du
	\right].
	\e*
Taking into account the terminal condition $\hat Y^\gamma_T = \partial_s g(\hat S^\gamma_T)$	and  the Euler Optimality condition for $\hat \alpha$ we get:~for all $\beta = (\beta^1, \cdots, \beta^\Gamma) \in \Ac^\Gamma$:
	\b*
	0&=&\Sum_\gamma \pi^\gamma  \E\left[  \int_0^T \left\{ \hat Y^\gamma_u +  
	P^{\bar {\hat \alpha}}_u + \partial_\alpha L^\gamma_T(Q^\gamma_u,\hat \alpha_u 
	+ \partial_\alpha L_S(\hat S_u, \hat \alpha_u)
	 \right.\right.\\&&
	 \hspace{20mm} \left. \left. -  p'\left(- Q^{0}_u - 
	 \Pi_{\Gamma} \cdot \left(\bar Q_u - \bar {\hat \alpha}_u \right)
	  \right)
	\left(  -  Q^{0}_u-\Pi_{\Gamma} \cdot \left(\bar Q_u - \bar{\hat  \alpha}_u \right)  \right)
	 \right\}\beta^\gamma_u  du\right].
	\e*
We deduce the coupling condition \reff{eq: coupling mfc}.

\begin{Proposition}
Assume that $\hat \alpha$ is a mean field optimal control for the  problem with a pricing rule $p$. Then $\hat \alpha$ is a mean field Nash equilibrium for the MFG problem with pricing rule
	\be
	p^{\rm MFG}(x)&=& p(x) + x p'(x)\;.
	\ee
\end{Proposition}
\proof
We first remark that the two McKean-Vlasov BSDEs   \eqref{eq: bsde} and \eqref{eq: bsde mfc}   are of the same form. Therefore  we obtain the desired result by comparing the two coupling equations \eqref{eq: coupling cnd}  and  \eqref{eq: coupling mfc}.
\ep

\section{The  Linear quadratic case}
\label{sect: explicit linear quadratic case}
In this section, we assume that the pricing rule is linear
	\be\label{eq: linear pricing rule}
	p:~x \mapsto p_0 + p_1 x.
	\ee
In this case, the function $\alpha \mapsto J^C_{x_0}(\alpha)$ is coercive and  strictly convex, which implies the existence of a unique mean field optimal control $\hat \alpha$.

\no Moreover, we assume that
	\b*
	\label{eq: quadratic storage cost}
	&&L_S:~ (s,\alpha) \mapsto \frac{A_2}{2} s^2 + A_1 s + \frac{C}{2} \alpha^2\\
	\label{eq: quadratic transmission cost}
	&&L^\gamma_T:~ (q,\alpha) \mapsto \frac{K^\gamma}{2}  \left( q - \alpha\right)^2\\
	\label{eq: quadratic terminal cost}
	&&g:~s \mapsto \frac{B_2}{2}\left(s-\frac{B_1}{B_2}\right)^2,
	\e*
where $p_0$, $p_1$, $A_1$, $A_2$, $C$, $B_1$, $B_2$ and $\{K^\gamma\}_{\gamma=1}^\Gamma$ are some given constants with $p_1>0$, $A_2>0$, $A_1<0$, $C<0$, $B_2>0$ and $K^\gamma\ge0  \quad \forall \gamma$.

\begin{itemize}
	\item In the storage cost $L_S$:~the term  $ (C/2)\alpha^2 $ is the current usage cost of the battery, it penalizes  the injection and withdrawal rate,  the term $(A_2/2) s^2$  is the current  cost of storage capacity and $A_1<0$ is the penalized negative stock level.
	\item  The demand charge $L^\gamma$ should be linked to the maximum instantaneous power consumed and is approximated in this setting by a quadratic expression.
	\item The terminal cost $g(S^{i, \alpha^i}_T)$ typically guarantees a a minimal level of storage at the end of the period.
\end{itemize}

{
In this setting, the convergence of the Nash-equilibrium for the $N$-player game to the EMFG can be proven. We mention that Graber \cite{Graber2016} (Section 3, Theorem 3.7,  p.15) shows for a class of linear-quadratic extended Mean Fields an approximate Nash equilibria property. Same arguments apply in our case and lead to the following convergence result.}
\begin{Proposition}[$\eps$-Nash equilibrium for the $N$-players game]
	Let  $\alpha^{i,\star}$ is a mean field Nash equilibrium for  $J^{\rm MFG}_{x^I_0}$. Then for each $\eps>0$ there exists $N_\eps$ and $\eta_\eps$ such that:  if  $N \ge N_\eps$ and $\eta \le \eta_\eps$,  then $\alpha^\star:=(\alpha^{1,\star}, \cdots, \alpha^{N,\star})$ is an $\eps$-Nash equilibrium for the N-players game.
\end{Proposition}

\subsection{Explicit solution of the MFC} $ $ 
Explicit solution of the MFC can be calculated. To simplify the notations in the following, we would write the optimal control of the MFC as simply $\alpha$ and not $\hat{\alpha}$ (the same for the controlled variable $\hat{S}$). To calculate the explicit solution, we first take the conditional expectation with respect to the common noise of the FBSDE system to calculate the optimal expected control $\bar{\alpha}$. Given this optimal expected control, the optimal control $\alpha$ can be caracterised in a second step.

Let's denote by $\hat{K}^\gamma:=C+K^\gamma$.  $\hat{K}^\gamma$ is strictly positive since we assume $C>0$ and $K^\gamma\ge0$.

Let also define the matrix 
$M_{\text{MFC}}:=\begin{pmatrix}
\hat{K}^1+2p_1 \pi_1 & 2p_1 \pi_2 & \cdots & 2p_1 \pi_\Gamma \\
2p_1 \pi_1 & \hat{K}^2+2p_1 \pi_2 & \cdots & 2p_1 \pi_\Gamma \\
\vdots&	\ddots& &\vdots\\
2p_1 \pi_1 & 2p_1 \pi_2 & \cdots &\hat{K}^\Gamma+2p_1 \pi_\Gamma \\

\end{pmatrix}. $ \\

Its determinant is $\text{det}_{{\text{MFC}}}=\prod_{j=1}^{\Gamma}(\hat{K}^j) +\sum_{j=1}^{\Gamma}(2p_1 \pi_j)\prod_{i\ne j}(\hat{K}^i).$  And its inverse matrix is $-M:=M_{\text{MFC}}^{-1}=\frac{1}{det_{M_{\text{MFC}}}} {\hat M_{\text{MFC}}}$ with $ {\hat M_{\text{MFC}}}$ the following matrix\\
\b*
 \begin{pmatrix}
	\Prod_{j\ne 1}\hat{K}^j + \sum_{j \ne 1}2p_1 \pi_j\Prod_{i\ne 1,j}\hat{K}^i& - 2p_1 \pi_2\Prod_{j \ne 1,2}\hat{K}^j &\cdots& -2p_1 \pi_\Gamma\Prod_{j \ne 1, \Gamma}\hat{K}^j\\

	-2p_1 \pi_1\Prod_{j \ne 1,2}\hat{K}^j & \Prod_{j\ne 2}\hat{K}^j + \sum_{j\ne 2}2p_1 \pi_j\Prod_{i\ne 2,j}\hat{K}^i&  \cdots& 2p_1 -\pi_\Gamma\Prod_{j\ne 2,\Gamma}\hat{K}^j\\
	
	\vdots&	\ddots& &\vdots\\
	-2p_1 \pi_1\Prod_{j \ne 1,\Gamma}\hat{K}^j&\cdots& &\Prod_{j\ne \Gamma}\hat{K}^j + \sum_{j \ne \Gamma}2p_1 \pi_j\Prod_{i\ne \Gamma,j}\hat{K}^i
\end{pmatrix} .
\e*
\bigskip

\no {\bf Step 1.}~
 In this linear quadratic case, if $\alpha$ is an optimal coordinated plan, we deduce from the  FBSDE \reff{eq: fsde}-\reff{eq: bsde mfc} and the  coupling condition \reff{eq: coupling mfc} that 
	\b*
	d\bar S_t&=& \bar \alpha_t dt,~~\bar S_0 = 0, \\ 
	d \bar Y_t&=& - (A_2\bar S_t  + A_1\1_\Gamma)dt + \bar Z^{0}_t dB^0_t,~~\bar Y_T = B_2 \bar S_T - B_1, \1_\Gamma\\
	\mbox{with}&&\bar \alpha_t\;= M\left(
	\bar Y_t +b_t	\right),
	\e*
where 
	\b*
	b_t = - \left( \diag{[\hat K_\Gamma ]} + 2p_1 \Pi_\Gamma \right) \bar Q_t -2p_1  Q_t^{0} \1_\Gamma + p_0 \1_\Gamma,
	\e*
 $\hat K_\Gamma=(\hat K^1, \cdots, \hat K^\Gamma)^T$ and  $M=- M_{\text{MFC}}^{-1}$.\\

By looking at a solution of the form $\bar Y_t - B_2\bar S_t = \bar \phi(t) \bar S_t + \bar \psi_t$, we are held to solve the following system:
\b*
\dot {\bar \phi} (t) + \bar \phi(t)M\bar \phi(t) + B_2 M \bar \phi (t) + B_2 \bar \phi(t) M + A_2 +B_2^2M &=& 0, \\
  \bar \phi(T)&=&0\\
d\bar \psi_t + (B_2 M + \bar \phi()t) M) \bar\psi_t dt +(\bar \phi(t) M b_t + B_2 M b_t + A_1 \1 _\Gamma)dt - Z_t^0 dB_t^0 &=&0, \\
 \bar \psi_T&=&-B_1 \1_\Gamma.
\e*

Denote by $\Ac = \left[ \begin{array}{cc}
B_2 M & M \\
-A_2 -B_2^2 M & -B_2 M
\end{array}\right]$
and referring to Theorem 5.3 in \cite{yong99} , if
\b*
\det \left[\left(0, I_\Gamma \right) e^{\Ac (T-t)}
\left(\begin{array}{c}
	0\\
	I_\Gamma
\end{array} \right)\right] >0
\e*
then $\bar \phi$ admits an explicit solution given by
\be
\bar \phi(t) &=& -\left[\left(0, I_\Gamma \right) e^{\Ac (T-t)}
\left(\begin{array}{c}
	0\\
	I_\Gamma
\end{array} \right)\right]^{-1}
\left[\left(0, I_\Gamma \right) e^{\Ac (T-t)}
\left(\begin{array}{c}
	I_\Gamma\\
	0
\end{array} \right)\right] .
\ee
 
 
By denoting $\chi_t$ the solution of the following linear ordinary differential equation
\b*
d\chi_t = (B_2 M + \bar \phi (t) M) \chi_t dt, ~~ \chi_0 = I_\Gamma, 
\e*
the solution of the linear BSDE for $\bar \psi$ is 
\begin{equation}
\bar \psi_t=  -\chi_t^{-1} \chi_T B_1 \1_\Gamma +\E\left[
\int_t^T  \chi_t^{-1} \chi_u \left((\bar \phi(uà M + B_2 M) b_u + A_1 \1_\Gamma \right)    du | \Fc_t^0 \right].
\end{equation}
 
Therefore, the optimal control $\bar \alpha$ is explicitly given because let's recall that $\bar \alpha_t\;= M\left(
\bar Y_t +b_t	\right)$ and $\bar Y_t - B_2\bar S_t = \bar \phi(t) \bar S_t + \bar \psi_t$. \\

\no {\bf Step 2.}~In this linear quadratic case, if $\alpha$ is an optimal coordinated plan, we deduce from the  FBSDE \reff{eq: fsde}-\reff{eq: bsde mfc} and the  coupling condition \reff{eq: coupling mfc} that 
\b*\label{eq: bsde mfc_step2}
d  S_t&=&  \alpha_t dt,~~ S_0 =  s_0,\\
d Y_t&=& - \left( A_2   S_t + A_1 \1_{\Gamma}\right) dt  +   Z^0_t dB^0_t +  Z_t dB_t, ~~ 	 Y_T= B_2 {S}_T-B_1 \1_\Gamma,\\
\mbox{with}&& {\alpha}_t = \hat M \left( Y_t
+\hat b_t\right), \\
&&   \hat M = \diag\left(\frac{-1}{C+K_\Gamma}\right) \\
\mbox{and}&& \hat b_t = p_0 \1_\Gamma - 2p_1( Q_t0r + \Pi_\Gamma(\bar{Q}_t-\bar{\alpha}_t^0))\1_\Gamma - \diag(K_\Gamma) Q_t .
\e* 

By looking at a solution of the form $ Y_t - B_2 S_t =  \phi(t)   S_t +  \psi_t$, we are held to solve the following system:
\b*
\dot { \phi} (t) +   \phi(t)\hat M\phi(t) + B_2 \hat M   \phi (t) + B_2  \phi(t) \hat M + A_2 +B_2^2\hat M &=& 0, \\
 \phi(T)&=&0\\
d  \psi_t + (B_2 \hat M +   \phi(t) \hat M)  \psi_t dt +(  \phi(t) \hat M  \hat b_t + B_2 \hat M \hat b_t + A_1 \1 _\Gamma)dt -  Z_t^0 dB_t^0 -  Z_t dB_t &=&0, \\
\psi_T&=&-B_1.
\e*

As $\hat M$ is diagonal, the solution of the Ricatti equation is explicit and by standard computations we can get
	\be
	  \phi^\gamma(t)+B_2&=&- \frac{\rho^\gamma}{\Delta^\gamma}\; \frac{
		e^{- \rho^\gamma (T-t)}(- B_2 \Delta^\gamma + \rho^\gamma) - e^{ \rho^\gamma (T-t)}(B_2 \Delta^\gamma + \rho^\gamma)
	}
	{
		e^{- \rho^\gamma (T-t)}(- B_2 \Delta^\gamma + \rho^\gamma) + e^{ \rho^\gamma (T-t)}( B_2 \Delta^\gamma + \rho^\gamma)
	}, \;\\
	\nonumber \mbox{with}~~&&\rho^\gamma := \sqrt{A_2 \Delta^\gamma},\\
	\nonumber &&\Delta^\gamma := \frac{1}{C+K^\gamma}.
	\ee
Let's define $ \phi^{\gamma,B_2}(t):= \phi^\gamma(t)+B_2$, then the solution of the BSDE is given explicitly by 
\be
  \nonumber \psi_t^\gamma= -B_1 \exp\left\{ - \int_t^T \Delta^\gamma \left(  \phi^{\gamma,B_2}(u)\right) du\right\}-\\
 \E\left[
\int_t^T  \Delta^\gamma    \phi^{\gamma,B_2}(u)  \exp\left\{ - \int_t^u \Delta^\gamma   \phi^{\gamma,B_2}(s) ds\right\}  \left(\hat{b}_u - \frac{A_1}{\Delta^\gamma  \phi^{\gamma,B_2}(u)} \right) du | \Fc_t 
\right].
\ee

Therefore, the optimal control $\alpha$ is explicitly given by these two previous equations because let's recall that $ \alpha_t\;= \hat M\left(
 Y_t + \hat b_t	\right)$ and $ Y_t - B_2 S_t =  \phi(t)  S_t + \psi_t$. \\

\subsection{Explicit solution of the MFC with 1 region} 
The case when the system is composed of only one region is already included in the results of the previous section. Nevertheless, we provide the expression in the special case of 1 region as the expressions are  more simple and we can also in that case express the controlled variable $S$. The system now becomes ($\pi=1$):

\no {\bf Step 1.}~ In this first step, we use the forward backward system  \reff{eq: fsde}-\reff{eq: bsde mfc} and the coupling condition \eqref{eq: coupling mfc} in order to get the optimal control $\bar{\alpha}$ and the optimal trajectory $\bar{S}$ associated to one node in this region.
\b*
d\bar S_t&=& \bar \alpha_t dt,~~\bar S_0 = 0,\\
d \bar Y_t&=& - (A_2\bar S_t  + A_1)dt + \bar Z^{0}_t dB^0_t,~~\bar Y_T = B_2 \bar S_T - B_1.
\e*
Rewriting the coupling condition \eqref{eq: coupling mfc}  in the case of one region gives
\begin{equation*}
\bar{Y}_{t}-K(\bar{Q}_{t}-\bar{\alpha}_{t})+C\bar{\alpha}_{t}+\bar{P}_{t}-p'(-Q^{0}_{t}-\bar{Q}_{t}+\bar{\alpha}_{t})(-Q^{0}_{t}-\bar{Q}_{t}+\bar{\alpha}_{t})=0,
\end{equation*}
where $\bar{P}_{t}$ is given by 
\begin{equation*}
\bar{P}_{t} =p^{MFG}(-Q^{0}_{t}-\bar{Q}_{t}+\bar{\alpha}_{t})=p_0+2p_1(-Q^{0}_{t}-\bar{Q}_{t}+\bar{\alpha}_{t}).
\end{equation*}
So, we obtain 
\begin{align*}
\bar{\alpha}_{t}&=-\frac{1}{K+C+p_1}\Big( \bar{Y}_{t}+p_0-p_1Q^{0}_{t}-(K+p_1)\bar{Q}_{t}\Big)\\
&=-\Delta(\bar{Y}_{t}+b_t)
\end{align*}
where $b_t=p_0-p_1Q^{0}_{t}-\bar{Q}_{t}(K+p_1)$ and $\Delta= \displaystyle \frac{1}{K+C+p_1}.$


By looking at solution of the form $\bar Y_t = \bar  \phi(t) \bar S_t + \bar \Psi_t$, we are held to solve the following system: $\bar \phi$ is the unique solution to the Riccati equation	
	\begin{equation} 
	\label{eq:riccati}
	\dot {\bar \phi} - \Delta \, \bar \phi^2 + A_2 =0 \; \mbox{with} \;  \bar \phi(T) = B_2
	\end{equation}
and $\bar \Psi$ is the unique solution to the linear BSDE	
	\begin{equation} 
	 \label{eq:LBSDE}
	d \bar \Psi_t=\Delta \bar \phi(t)  \left(  \bar\Psi_t + \bar P_t  \right) dt + \bar Z^{0}_t dB^0_t,~~\bar \Psi_T = -B_1.
	\end{equation}
Using the affine form of the solution $ \bar Y$, the  Ricatti equation \eqref{eq:riccati} and by identification with equation \eqref{eq:LBSDE}, we obtain the explicit expression of the energy price 
 \begin{equation} 
	 \label{eq:price}\bar P_t = - \frac{A_1}{\Delta \bar \phi} + b_t.
	\end{equation}
Moreover, by using the standard computations, we get from \eqref{eq:riccati} and  \eqref{eq:LBSDE} the expressions of  $\bar \phi$  and  $\bar \Psi$ as following	\b*
	\bar \phi(t)&=&- \frac{\rho}{\Delta}\; \frac{
	e^{- \rho (T-t)}(- B_2 \Delta + \rho) - e^{ \rho (T-t)}(B_2 \Delta + \rho)
	}
	{
	e^{- \rho (T-t)}(- B_2 \Delta + \rho) + e^{ \rho (T-t)}( B_2 \Delta + \rho)
	}\;~~\mbox{with}~~\rho := \sqrt{A_2 \Delta},
	\e*	
	\begin{equation*}
	\bar \Psi_t= - B_1 \exp\left\{ - \int_t^T \Delta \bar \phi(u) du\right\}-\E\left[
	 \int_t^T  \Delta \bar \phi(u)  \exp\left\{ - \int_t^u \Delta \bar \phi(s) ds\right\} \bar P_u du | \Fc^0_t 
	 \right].
	\end{equation*}

It follows that $\bar S_t$ satisfies
\begin{equation*}
\bar S_t=
-\Delta \int_0^t    \exp\left\{ - \int_u^t \Delta \bar \phi(s) ds\right\} \left(
\bar P^{}_u + \bar \Psi_u + \frac{ A_1}{\ \Delta \bar \phi(u)}\right)
 du .
\end{equation*}
\no {\bf Step 2.}~Once we obtain all the optimal elements of one node in the first step, we use the FBSDE \reff{eq: fsde}-\reff{eq: q}-\reff{eq: bsde} and the coupling condition  \reff{eq: coupling cnd} to find the optimal objects associated to  one region containing a number of identical nodes. Thus, following the same computations  as the first step and looking for $ Y_t = \varphi(t) S_t + \psi_t$,   we get:
	\b*
	\alpha_t&=&  - \delta \left( Y_t + P_t + \frac{A_1}{\delta \phi(t)}.  \right),
	\e*	
where
	\b*
	\delta = \frac{1}{C+K}&\mbox{and}&
	P_t= p_0  - 2 p_1 (  Q^{0}_t + \bar Q_t - \bar \alpha_t) - K Q_t -\frac{A_1}{\delta \phi(t)}.
	\e*
Then the FBSDE \reff{eq: fsde}-\reff{eq: bsde} becomes
	\b*
	dS_t&=& - \delta \left( Y_t + P_t + \frac{A_1}{\delta \phi(t)}\right) dt, ~~S_0 = s_0,\\
	dY_t&=& - (A_2 S_t + A_1) dt + Z^0_t dB^0_t + Z_t dB_t,~~Y_T = B_2 S_T - B_1.
	\e*
Again, we  get  explicitly in the same ways as the first step
\b*
	\varphi(t)&=&- \frac{\rho}{\delta}\; \frac{
	e^{- \rho (T-t)}(- B_2 \delta + \rho) - e^{ \rho (T-t)}(B_2 \delta + \rho)
	}
	{
	e^{- \rho (T-t)}(- B_2 \delta + \rho) + e^{ \rho (T-t)}( B_2 \delta + \rho)
	}\;~~\mbox{with}~~\rho := \sqrt{A_2 \delta},
	\e*
	\begin{equation*}
	 \psi_t=- B_1 \exp\left\{ - \int_t^T \theta  \varphi(u) du\right\}-\E\left[
	 \int_t^T  \theta  \varphi(u)  \exp\left\{ - \int_t^u \theta  \varphi(s) ds\right\}  P_u    du | \Fc_t 
	 \right],
	\end{equation*}
and 
\begin{equation*}
S_t= s_0 \exp\left\{ - \int_0^t \delta  \varphi(u) du\right\} 
- \delta \int_0^t  \exp\left\{ - \int_u^t \delta \varphi(s) ds\right\} \left( 
P_u +  \psi_u +\frac{A_1}{\delta \phi(u)}
\right) du .
\end{equation*}	

\section{Numerical interpretations}
\label{sect:  numerical}

\subsection{Description of the game}
\label{subsect: description_game}

The "rest of the world" region is composed by agents who are traditional consumers and do not consider the opportunity to have storage and just face random consumption for electricity and pay the resulting random bill for their electricity. Indeed, their consumption is random but also spot prices they pay for their energy. The prosumer zones gather prosumers who optimize the capacity size of an individual battery and their injections and withdrawals. Agents could be consumers, producers or alternatively both. This last situation may represent residential consumers with photovoltaic pannels on top of their roof. These Agents are indeed producers during daytime when the sun shines and while they are out for work and these same Agents are consumers when they get back home at sunset. We will consider examples with one or two prosumer zones with different characteristics like their load demand volatility, seasonality... \\

{\bf Remark: possible extension of the proposed model to demand side management.}~Let's point out that our model can be extended to handle demand side management with respect to little adjustments. Indeed, demand response actions mainly consist in postponing or moving forward electricity usages that can be typically represented by a storage. Costs of storage then represent the costs of effort it takes to the Agent to modify its electricity load demand. \\

The optimization horizon $T$ of the Agents is typically several hours like a day or two. Indeed, we have in mind that residential batteries we represent in our problem can help to dispatch Agent's consumption over this horizon but not longer. In the simulations, we consider $T=1$ day.\\

{\bf Remark: model parameters.}~Our examples are designed to illustrate some stylized behaviors of the model and parameter values we use in the following are not based on real figures. Extensions of our model should be considered in the future, in particular the illustration of a real system.

The random injection of the prosumer and rest of the world zones are modeled as the sum of a deterministic seasonal function $\mu$ and an Ornstein-Uhlenbeck (OU) process (without independent noise for the rest of the world zone).
	\b*\label{eq: node i Q2}
	dQ^{i}_t&=& -a^{\gamma} (Q^{i}_t-\mu^{\gamma}(t)) dt  + \sigma^{\gamma} dB^{i}_t + \sigma^{\gamma~0} dB^0_t\;,~~Q^{i}_0 = q^{i}_0,~~ i\in\gamma,\\
	dQ^{0}_t&=& -a^{0} (Q^{0}_t-\mu^{0}(t)) dt  + \sigma^{0} dB^0_t\;,~~Q^{0}_0 = q^{0}_0.
	\e*
 We consider here only one prosumer zone, ie $\Gamma = 1$. We will consider examples in the following where the seasonality of the rest of the world is twice in average the one of the prosumer zones. The seasonality $\mu$ is a simple cosine function which is a proxy for the peak and off-peak consumption of residential Agents. To summarize, the seasonal component of the consumption are given for each date $t$ expressed in day by:
	\b*
	\mu^0(t) = 2\cos(4\pi t - \pi/2)-3 \text{ and } \mu^\gamma(t) = \mu^0(t)/2.
	\e*
The other parameters of the model, if not stated otherwise, are in the following of the analysis: $a^0=a^\gamma=1$, $\sigma^\gamma=\sigma^0=0.8$ ,  $\sigma^{\gamma,0}=0.3$, $p_0=5$, $p_1=5$, $A_2=250$, $A_1=-15$, $C=5$, $K=10$, $B_2=5000$ and $B_1=-0.12B_2$.\\
Next figure is an example of random trajectories of the consumptions of several Agents with corresponding spot prices driven by linear pricing rule \ref{eq: linear pricing rule}. It happens that consumption can be negative which means that the Agents are producing electricity at that particular time. In the meanwhile spot prices can be negative which is an observed feature of electricity spot price which typically occurs when the residual consumption (consumption minus wind/solar productions) is very low, see for example \cite{Paraschiv}.\\

\begin{figure}[!ht]
	\begin{center}
		\hspace*{-6mm}
		\includegraphics[scale=0.35]{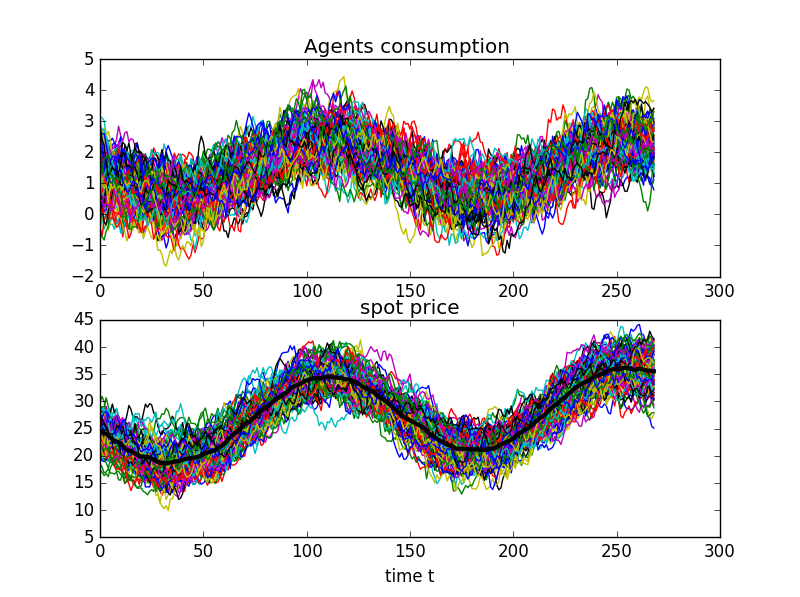} 	
		\caption{\label{fig: ex8Price} Agent's consumption (upper figure) and corresponding spot price (low figure) with average prices (wide black line) for several simulations, $T=1$ day.}
	\end{center}
\end{figure}

\subsection{Management of the storage with respect to the bill structure and the impact of Agents on spot price}

To have a storage enables Agents to influence two part of their electricity bill: 
\begin{itemize}
	\item {\bf to reduce the cost of the volumetric part} of their electricity bill by reporting their consumption/production when spot prices are low/high which also means time-arbitraging spot. By doing so, they have a smoothing impact on spot prices: their peak consumption is shifted during low global demand period whereas their off-peak demand is shifted during high global consumption period. This impacts directly the other population who also pays their volumetric part at the spot price.
	\item {\bf to reduce the cost of their capacity charge} by limiting their maximum load demand. In general, this has less influence on the other consumer, ie. the rest of the world, as this smoothes less spot prices.\\
\end{itemize}

Several factors imply that Agents are going to use their storage rather to favor one reduction or the other:
\begin{itemize}
	\item {\bf The influence of the Agent's consumption on the spot price}: the influence of the Agent's consumption is measured by two factors. First of all is the individual impact of the Agents linked to the size of the region he belongs to with respect to others (represented by parameter $\pi_i$). The second factor is the price differential between peak and off-peak period linked in our model to the global influence of the electricity consumption of the whole system over the spot price represented by parameter $p_1$. Small dissemination of storage in the system (low $\pi_i$) and/or large peak/off-peak spot price differential (high $p_1$) favor spot arbitrage and the willingness by the Agent to use their storage to reduce the cost of their volumetric part. Indeed, high $\pi_i$ which means lots of storage on the sytem  will diminish the interest of storage to make spot arbitrage because for example the individual Agent who decides to store to benefit from low spot price is also imitated by many others which has for consequence to increase spot price. On the contrary, low $p_1$ implies that the seasonality of spot price is less and automatically reduces the peak/off-peak differential. 
	
	\item {\bf The bill structure}: depending on the proportional weight of the volumetric part of the bill ($P^{N,\alpha}_t \left(Q^i_t - \alpha^i_t\right)$) compared to the demand charge part of the bill ($\frac{K^\gamma}{2}|Q^i_t - \alpha^i_t|^2 $), the Agents manage their storage differently. If bills are driven mainly by the demand charge, ie high $K$, the Agents use their storage so that they smooth the seasonality of their consumption and even obtain a residual load $-Q^i_t + \alpha^i_t$ nearly constant and as close as possible to  the average of the load $-Q^i_t$ over the period.	
\end{itemize}

\begin{figure}[!ht]
	\begin{center}
		\hspace*{-6mm}
		\includegraphics[scale=0.35]{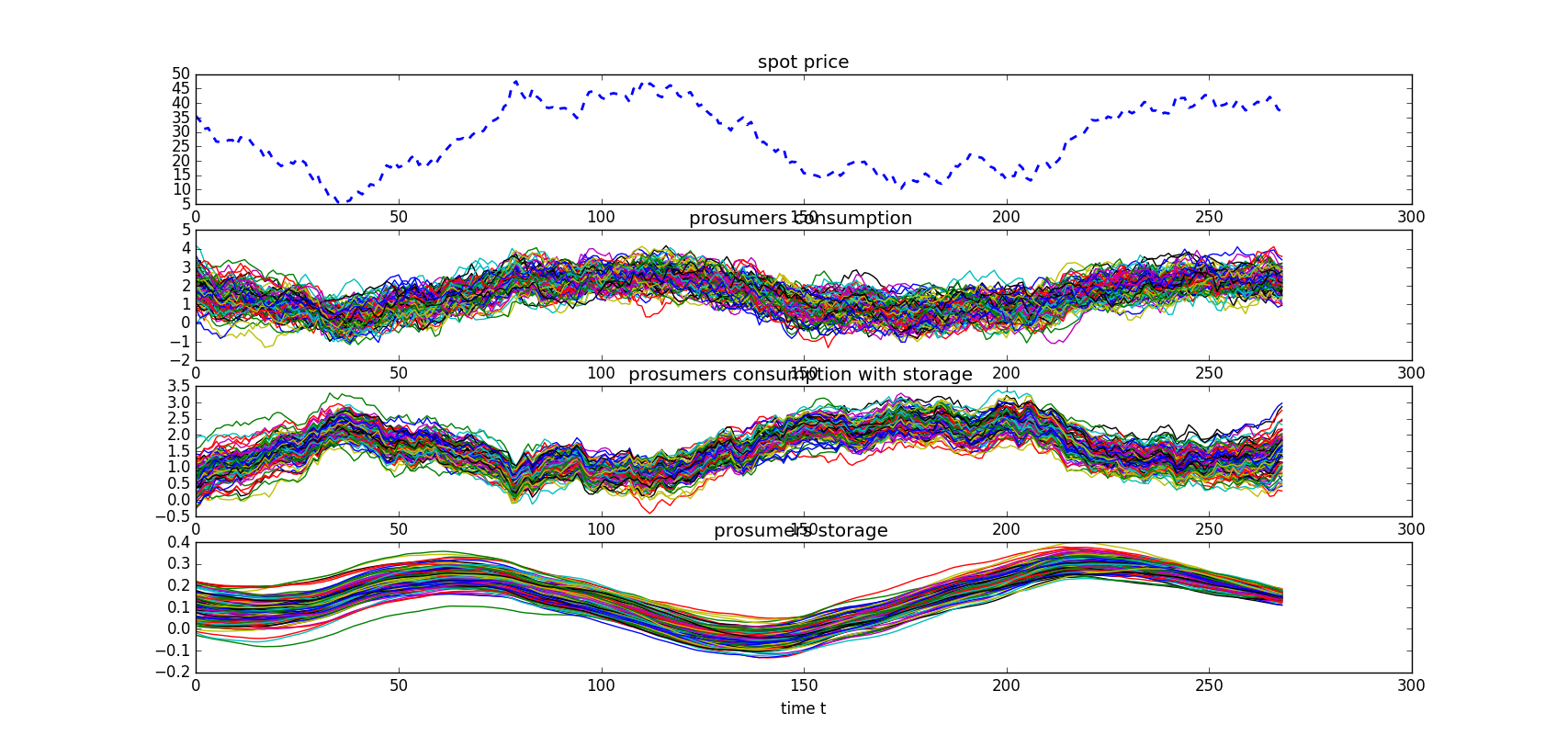} 	
		\caption{\label{fig: ex8Stock} One simulation of spot price (upper graph), prosumers' consumption $Q^i$ (middle graph), prosumer's net consumption $Q^i-\alpha^i$ (lower middle graph) and prosumer' storage level (lower graph) for every prosumers.}
	\end{center}
\end{figure}

Let's illustrate these conclusions by numerical examples. First, we modeled the rest of the world and the prosumer zones to be equivalent in terms of consumption but we suppose that the {\bf prosumers' zone has no  influence over the spot price} compared to the traditional consumers zone. It means that even if the number of prosumers is non negligible (can even be approximated as being infinite), their number compared to traditional consumers is low. This should correspond to a situation where residential storages have being developed but are still an exception in the population. Fig. \ref{fig: ex8Price} shows one simulation of spot price and the consumption $-Q^i$ of the prosumers before they consider using their storage. The optimized way to used their storage is, as expected, to store when prices are low and to withdraw when prices are high as shown in Fig. \ref{fig: ex8Stock} on one simulation of spot and Agents' consumption.\\

Let's point out that the storage curves are almost always positive. Negative values do occur but do not disrupt interpretations we can deduce from the model. Indeed, these negative value may be thought as the necessity to consider an energy reserve in the storage: the storage in normal mode is always operated above an energy reserve which may be necessary to use for some particular consumption/storage level occurrences. 

As expected, the resulting consumption that prosumers are adressing to the network is therefore a mirror of their initial ones as shown in Fig. \ref{fig: ex8Stock} (compare the two middle graphics). The storage is used such that Agents are reporting their high consumption when prices are low and are consuming less when prices are high. In addition, their net load demand $-Q^i + \alpha^i$ is smoother compared to original consumption $-Q^i$. To have local storage enable to reduce the maximum instantaneous power consumption in average by 21\% for every prosumers and reduces the electricity bill of prosumers by more than 13\% (the total reduction after including storage costs is only 7\%). This is summarized in the following table which indicates the repartition of the bill between the volumetric part and the demand charge part and the reduction on both parts implied by having a local storage.\\

\begin{equation*}
\begin{array}{|c|c|c|}
\hline
	& \text{electricity bill } & \text{reduction implied by battery} \\
	\hline
	\text{volumetric charge} & 76 \% & 21\% \\
	\text{demand charge} & 24 \% & 8\% \\
	\hline
\end{array}
\end{equation*}
\begin{center}
	prosumers - battery owners
\end{center}

\begin{figure}[!ht]
	\begin{center}
		\hspace*{-6mm}
		\includegraphics[scale=0.35]{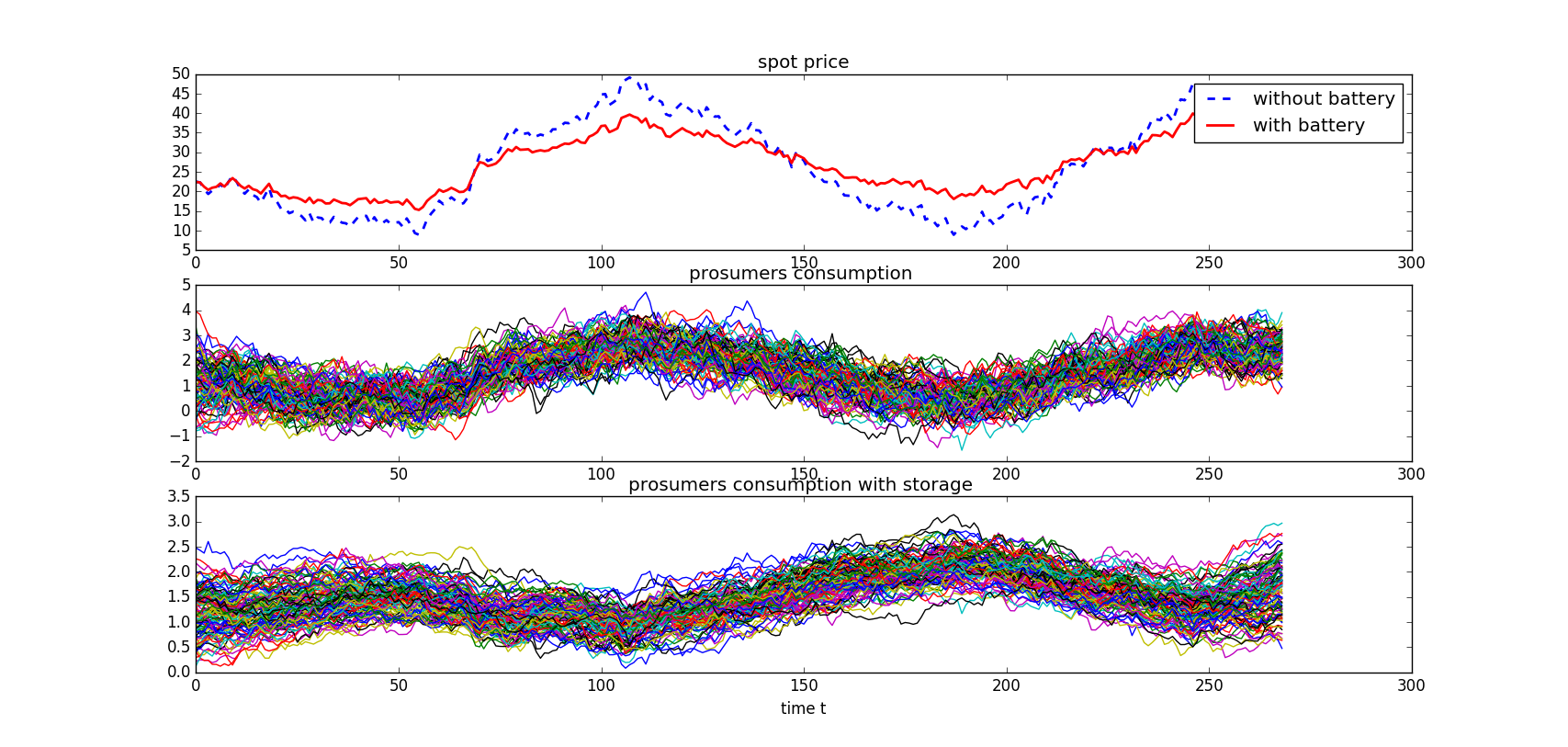} 	
		\caption{\label{fig: ex8ConsoInfluence} One simulation of spot price (upper graph) without battery in the system (straight line) and with batteries (dashed line), prosumers' original consumption $Q^i$ (middle graph), and prosumers' net consumption  $Q^i-\alpha^i$ (lower graph) for every prosumers.}
	\end{center}
\end{figure}

Let's now study the case when the {\bf prosumer zone has now equal influence on the spot price} as the rest of the world, which means that batteries would have spread among the population in such a way that battery owners and non-battery owners are equally distributed among the total population. In this case, the benefit of having a local battery is as expected slightly lower. Indeed, to postpone a large consumption when price are lower is less efficient because every prosumers do the same and as such make the spot price increase. We now observe the following impacts (after having modified spot price parameter $p_1$ such that the average spot price  remains the same as the previous one).
\begin{itemize}
	\item The spot price are smoothed (maximum prices decrease whereas minimum prices increase) and their volatility decreases (see upper graph of fig. \ref{fig: ex8ConsoInfluence}). This smoothing benefits to non-storer zone, indeed the spot price diminishes when their consumption is high and spot price increases when their consumption is low which has a lower impact on their bill. The "rest of the world" bill has diminished by 5 \%.
	\item It is not optimal, contrary to previous example, to completely flip the maximum and minimum consumption using the battery (see middle and lower graph of fig. \ref{fig: ex8ConsoInfluence}), as such the reduction of the electricity bill on the volumetric charge is lower than in the previous case when the influence on the spot price of prosumers was very low,
	\item The prosumers make more effort to gain on their demand charge part of their bill: they diminish their maximum consumption more (30\% reduction compared to 21 \% reduction when they have no influence on the spot price) because their main interest is no more spot arbitrage.  The impact on the prosumer's bill is given in the following table.
	\item The optimal battery capacity is slightly lower. \\
\end{itemize}

\begin{equation*}
\label{prosumers - battery owners}
\begin{array}{|c|c|c|}
\hline
& \text{electricity bill } & \text{reduction implied by battery} \\
\hline
\text{volumetric charge} & 76 \% & 13\% \\
\text{demand charge} & 24 \% & 16\% \\
\hline
\end{array}
\end{equation*}
\begin{center}
	Impact on electricity bill for battery owners
\end{center}

{\bf Remark for autosufficient prosumer:}~we observe that a prosumer who produces in average enough to fulfill its consumption in energy can disconnect from the system if the gain on spot is too little. \\

\subsection{Impact of decentralized management of batteries against centralized management}
The impact of decentralisation against centralization optimization can be measured with the common notion in game theory of Price of Anarchy, PoA. PoA measures the ratio of the total costs of all zones obtained with decentralized optimization (MFG optimization) an the costs of the total costs of all zones obtained with centralized optimization (MFC optimization)). PoA is always greater than 1. \\

\begin{figure}[!ht]
	\begin{center}
		\hspace*{-6mm}
		\includegraphics[scale=0.25]{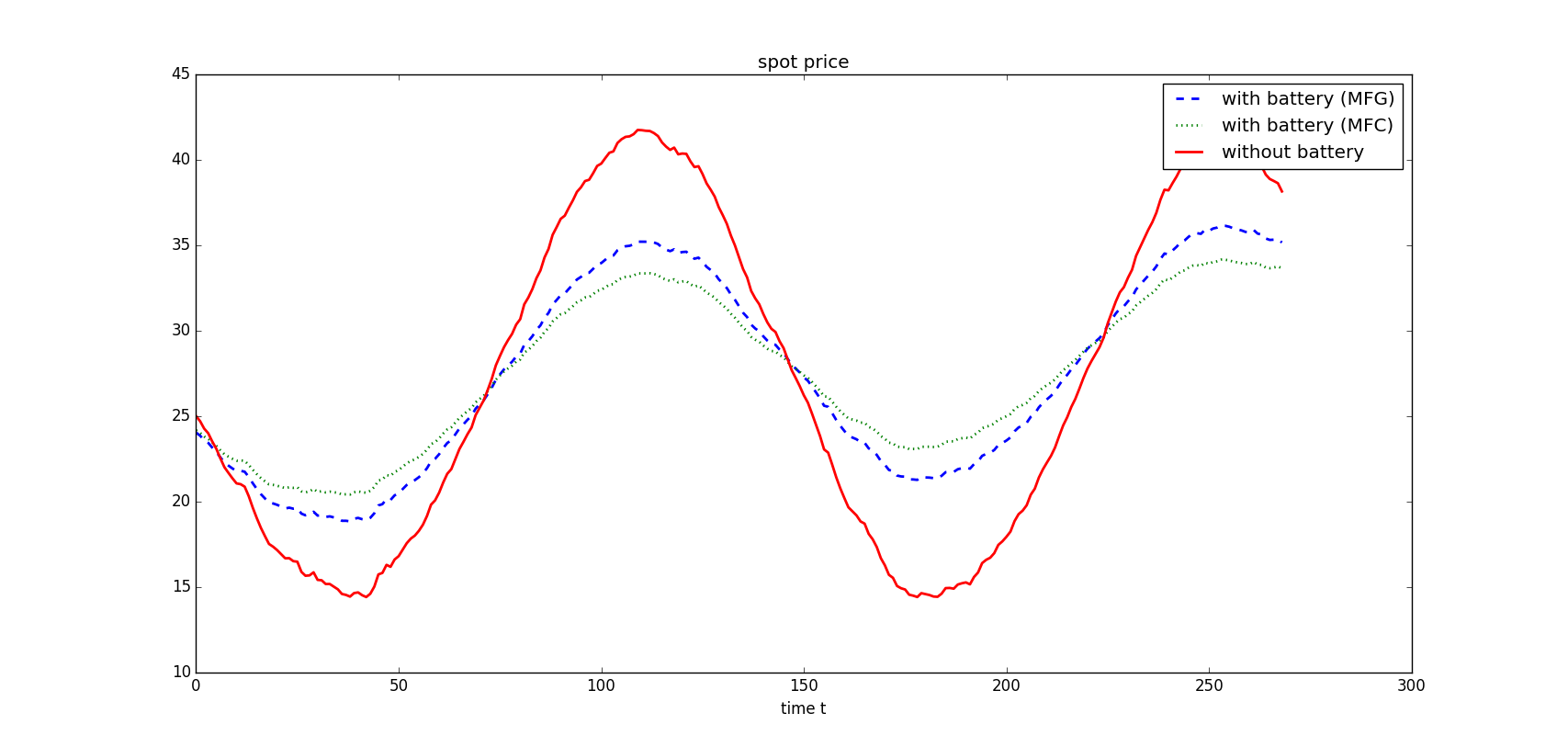} 	
		\caption{\label{fig: ex9Price} Average spot price over without battery in the system (straigth line), with decentralized batteries (dashed line) and with batteries optimized by a central planner (dotted line)}
	\end{center}
\end{figure}

In the example  we consider,  with two equivalent zones in terms of consumption and influence on the spot price, PoA is close to 1 meaning that the impact of having decentralized batteries in the system for the two consumer zones is not too high and that the optimization is rather close to what would be obtained by a centralized planner. Nevertheless, we observe some slight impacts: indeed a centralized management would allocate cost reductions more in favor to normal consumers ("rest of the world") than what a decentralized management does.
\begin{itemize}
	\item A centralized planner would install slightly higher battery capacity which would penalized a bit the battery owners zone because the cost of their battery would increase.
	\item To have bigger batteries would make the spot price smoother (see fig. \ref{fig: ex9Price}) and benefit  the population without battery by reducing more their energy payment (7\% cost reduction compared to 5\% in decentralised MFG management).
\end{itemize}

\subsection{Load demand variability increases the benefit of storage}

The more the volatility of the load demand, the more useful the batteries are for prosumers. Indeed, when the volatility of load increases, the fraction of the bill related to the demand charge increases. If the consumption variability is 2.5 times higher, the battery still diminish the maximum consumption power by around 30 \%, this has therefore a bigger impact on the bill (22\% reduction to be compared to 13 \% when standard consumption variability). Of course, in order to be able to reduce the maximum capacity of the prosumer's consumption in the same order as when the volatility of its consumptions is 2.5 times lower, the battery capacity also increases with the variability of consumption. To summarize,  increase of load variability has two main impacts:
\begin{itemize}
	\item increase of battery capacity of prosumers,
	\item a larger reduction of the electricity bill. 
\end{itemize}

\begin{equation*}
\begin{array}{|c|c|c|}
\hline
& \text{electricity bill } & \text{reduction implied by battery} \\
\hline
\text{volumetric charge} & 67 \% &  15\% \\
\text{demand charge} & 23 \% & 28\% \\
\hline
\end{array}
\end{equation*}
\begin{center}
	Impact on electricity bill for battery owners and system with 2.5 higher consumption volatilities
\end{center}

\subsection{Example of two prosumer competing zones}
Our model can deal with several prosumers' zones. Let's modify a bit our core example to illustrate a competition between two zones. We consider now one prosumer zone whose seasonal pattern of consumption is in opposition with the "rest of the world" . This means that the prosumer peak consumption now occurs when the "rest of the world" has its lowest consumption. Without storage, spot price pattern is still governed by the "rest of the world" consumption seasonality (because the seasonality of "rest of the game" is twice the one of prosumer zone as chosen in section \ref{subsect: description_game}). This induced that the energy cost of prosumers, without storage, is now lower than in previous examples (only 70\%) because they naturally consumes when prices are the lowest. \\

If this prosumer zone now installs local batteries, prosumers will install lower battery capacity than in previous examples and only  fulfill the objective to diminish their demand charge (indeed their consumption pattern is naturally optimal and their benefit from spot arbitrage is then very low). By doing so, prosumers reduce their maximum consumption which occurs at off-peak and therefore reduce the off-peak spot price slightly. This reduction of consumption is reported when their consumption is at the lowest which also corresponds to the peak of spot prices and therefore  makes the peak spot price slightly increase. In this example, the storage management has a negative impact for the "rest of the world population" which has its energy part of its bill slightly increases (1\% increase). \\

We show by numerical simulation that if the prosumer zone is now divided in two zones, $\Gamma=2$, of equal size: one zone with a seasonal pattern in phase with the "rest of the world" and referred next as "in-phase" zone (studied in previous subsections) and one in opposition to the seasonality of the "rest of the world" and referred as "de-phase" zone (studied in above in this subsection). In that case, the "de-phase" zone will suffer an increase of its bill after having installed batteries because the "in-phase" has also installed batteries. By doing so, the "in-phase" zone has smooth spot prices which is negative for "de-phase" zone. The "de-phase" zone would then lose from battery installation in the system (whereas it is still beneficial to install batteries for the "de-phase" zone or it would loose even more).

\subsection{Conclusion of numerical tests}
Examples presented in this paper are some illustrations of what the model can enable to study. Many other experiments and tests can be conducted easily because the model is quite generic. Let's recall that the model can cope with quite general dynamics for the consumptions/production and is not limited to the simple Ornstein-Uhlenbeck considered here. In particular, the implementation of cases calibrated on real figures should be conducted in future research. Very recently \cite{Huetall17} caracterised MFG with constrained controls, their results may be applied for our class of Extended-MFG to study how physical constraints of the storage influence numerical results.

\paragraph{Acknowledgements} The authors wish to thank the anonymous referee for all the pertinent remarks she/he made.
The authors's research is part of the ANR project CAESARS (ANR-15-CE05-0024) and PACMAN (ANR-16-CE05-0027) and of PANORISK project. 
The third author was partially supported by chaire Risques Financiers de la fondation du risque, CMAP-Ecole Polytechniques, Palaiseau-France.

\bibliographystyle{siam} 
\bibliography{Extended-MFG-18} 
\end{document}